\documentclass[12pt]{article}
\usepackage{amsmath,amssymb,amsfonts,amsthm}
\usepackage[cp1251]{inputenc}
\usepackage[T2A]{fontenc}
\usepackage[english,ukrainian]{babel}
\newtheorem{theorem}{Theorem}

\begin{document}
\def\refname{References}
{\bf REGULARIZATION OF LINEAR IMPULSIVE BOUNDARY-VALUE PROBLEMS FOR SYSTEMS OF
INTEGRO-DIFFERENTIAL EQUATIONS}\\

\begin{center}
Ivanna BONDAR\\
Institute of Mathematics of NASU, Kyiv, Ukraine\\
\textit{bondar.ivnn@gmail.com} 
\end{center}

\textbf{Abstract.}{\small \it \noindent The solvability of impulse system of integro-differential equations with a degenerate kernel is investigated.
It is assumed that the impulse system does not have a solution for arbitrary inhomogeneities.
In order to reduce it to solvable, a control function was introduced, a solvability criterion
was established, and its general form was constructed. The fact that the control may not be unique allows us
to use it to study problems that are often encountered in the theory of optimal control.
The general method of studying the problem posed in this way uses the theory of pseudoinverse matrices (in the Moore–Penrose sense) and orthoprojectors.
\hfill}\\



\begin{flushleft}
\textbf{1. Introduction}
\end{flushleft}
In various applied sciences, mathematical models of processes are appearing, which are described by systems of
algebraic and integro-differential equations (IDE). A wide range of such mathematical models are described by systems of
IDE with various kinds of disturbances or nonlinearities. It is known that some problems of
optimal control, linear programming, economics, theory of elasticity, hydrodynamics, chemical and biological kinetics, etc.
are modeled by such operator equations. In researching the solvability of various types of functional differential equations
and boundary value problems for them, the theory of generalized inverse operators \cite{BoiSam, BLS,ZhurSl} has been widely used in the last decade.
This approach allows, taking into account the specifics of each specific problem, to apply all the advantages of the "operator theory"
for its solution. The specificity of studying the solvability and construction of solutions of IDE is that their linear part is an operator that does not have an inverse.

This fact significantly complicates the study of such operator equations and BVPs for them and leads to the fact that the solution of the BVP for such systems consists of the solvability conditions of both the operator system itself and the boundary value problem for it \cite{BoiHol,LanMar,SBK}.

To investigate the existence of solutions to such problems, as will be shown below, you can use the apparatus of the theory of pseudo-inverse matrices
and operators, which was developed in the works of A.M. Samoilenko, O.A. Boichuk \cite{BoiSam,SamPer} and actively developed for the case of weakly perturbed boundary value problems for systems of integro-differential equations \cite{Bon_7}, impulse boundary value problems \cite{BonGrKoz, Bon_16} and BVPs for
integro-dynamic equations on time scales \cite{BonNestStr}.

\begin{flushleft}
\textbf{2. Main result}
\end{flushleft}
Studying control conditions for impulse boundary value problems (IBVP) is crucial when dealing with systems of integro-differential equations (IDE). These conditions play a fundamental role in determining the behavior and stability of the system, as well as the existence and uniqueness of solutions. The basis of this study was work \cite{BoiHol}, which investigated solvability of linear system of IDE with a degenerate kernel.
That is why we consider IBVP for IDE and assum that it has no solution for arbitrary inhomogeneities. In order to reduce it to a solvable one, we are going to introduce a control function, establishe a solvability criterion and construct a general form for a solution of it.

Consider the inhomogeneous system of integro-differential equations with impulsive actions at fixed times
\begin{equation}\label{b:1}
\begin{array}{c}
  \dot{x}(t)-\Phi(t)\int\limits_{a}^b
\Big[A(s)x(s)+B(s)\dot{x}(s)\Big]ds=f(t)+\int\limits_{a}^b
K(t,s)ds\cdot u,~~t\neq\tau_{i},\\
 \Delta E_{i}x\mid_{t=\tau_{i}}:=S_{i}x(\tau_{i}-0)+\gamma_{i},~~i=1,...,p,
\end{array}\end{equation}

\begin{equation}\label{b:2}
\ell
x(\cdot,\varepsilon)=\alpha~\in~\mathbb{R}^{q}.
\end{equation}

Here, we use the assumptions and notation from
\cite{BoiSam,BoiHol}: $A(t),~B(t),$ $\Phi(t),$ $K(t,s),K_1(t,s)$
are, respectively, $m\times n,$ $m\times n,$ $n\times m,$ $n\times
n,$ $n\times n$ matrices~ which components are sought in the space
$L_2[a,b];$ column vectors of matrice $\Phi(t)$ are linearly
independent at $[a,b];$ the $n \times 1$  vector function $f(t) \in
L_{2}[a,b]$; $E_{i},~S_{i},~A_{1i}$ are $k_{i}\times n$ constant
matrices such that rank$(E_{i}+S_{i})=k_{i}<n,$ which means that the
corresponding components of solutions of the impulsive system admit
unambiguous continuation through the points of discontinuity
$$\Delta E_i x\Bigr|_{t=\tau_i}:= E_i ( x(\tau_{i}+
0)-x(\tau_{i}-0));$$ $\gamma_{i}$ is an $k_{i}$-dimensional column
vector of constants, $\gamma_{i} \in \mathbb{R}^{k_{i}};$
$$a<\tau_{1}<...<\tau_{i}...<\tau_{p}<b~~~\texttt{for}~~~i=1,...,p;$$
$\ell=col(\ell_{1},\ell_{2},...,\ell_{q})$ is a bounded linear
$q$-dimensional vector functional,\\
$\alpha={\rm col}(\alpha_{1},\alpha_{2},...,\alpha_{q}) \in
\mathbb{R}^{q}.$

The solution $x(t)$ is sought in the space of $n$-dimensional
absolutely continuously differentiable vector functions
\begin{center}
$x=x(t,\varepsilon):x(\cdot,\varepsilon)\in
D_{2}([a,b]\setminus_{\{\tau_{i}\}_{I}})$,$~~~\dot{x}(\cdot,\varepsilon)\in
L_2[a,b]$,$~~~ x(t,\cdot)\in C(0,\varepsilon_0].$
\end{center}
The norms in the spaces $D_{2}([a,b]\setminus_{\{\tau_{i}\}_{I}})$
and $L_2[a,b],~C(0,\varepsilon_0]$ are introduced in the standart
way (by analogy with \cite{BoiSam, ZhurSl}).

We first consider that problem (\ref{b:1}), (\ref{b:2}) belongs to
the Fredholm case ($n\neq p$) and obtain bifurcation conditions of
solution of  this problem from the point $\varepsilon = 0$.

Parallel with the weakly perturbed  boundary-value problem
(\ref{b:1}), (\ref{b:2}), we consider the following generating
boundary-value problem ($\varepsilon=0$):

\begin{equation}\label{b:3}
\begin{array}{c}
  \dot{x}(t)-\Phi(t)\int\limits_{a}^b
\Big[A(s)x(s)+B(s)\dot{x}(s)\Big]ds=f(t),~~t\neq\tau_{i},\\
 \Delta E_{i}x\mid_{t=\tau_{i}}:=S_{i}x(\tau_{i}-0)+\gamma_{i},~~i=1,...,p,
\end{array}\end{equation}

\begin{equation}\label{b:4}
\ell x(\cdot,\varepsilon)=\alpha~\in~\mathbb{R}^{q}.
\end{equation}
Assume that the boundary-value problem (\ref{b:3}), (\ref{b:4}) does
not  have solutions for  arbitrary  inhomogeneities $f(t)\in L_2
[a,b]$ and $\alpha \in \mathbb{R}^{p}$.

Impulsive condition can be written as the interface boundary
conditions \cite{Zet} by using the $k$-dimensional linear bounded
vector functional
\begin{equation*}\begin{array}{c}
\varphi=col(\varphi_{1},\varphi_{2},...,\varphi_{p}):~D_{2}([a,b]
\setminus_{\{\tau_{i}\}_{I}})\rightarrow
\mathbb{R}^{k},\\
\varphi_{i}:~D_{2}([a,b]\setminus_{\{\tau_{i}\}_{I}})\rightarrow
\mathbb{R}^{k_{i}},\\
k:=k_{1}+k_{2}+...+k_{p},~~i=1,2,...,p
                 \end{array}\end{equation*}
where
\begin{equation} \label{b:5}
\left\{\begin{array}{l} {\varphi_{1} x:=E_{1}x(\tau _{1} +)-(E_{1}+S_{1} )x(\tau _{1} -)} \\
{\varphi_{2} x:=E_{2}x(\tau _{2} +)-(E_{2}+S_{2} )x(\tau _{2} -)} \\
{................................................} \\
{\varphi_{p} x:=E_{p} x(\tau _{p} +)-(E_{p}+S_{p} )x(\tau _{p} -)}
\end{array}\right.
\end{equation}
and have next form
\begin{equation}\label{b:6}
\varphi x(\cdot,\varepsilon)= \gamma \in \mathbb{R}^{k}.
\end{equation}
Here $\gamma={\rm col}(\gamma_1, \gamma_2, ..., \gamma_{p}) \in
\mathbb{R}^{k},~~\gamma_{i} \in \mathbb{R}^{k_i}.$

We introduce the bounded linear $(k+q)$-dimensional vector
functional
$\mathfrak{L}:=[\begin{array}{c}
                  \varphi \\
                  \ell
                \end{array}]:D_{2}([a,b]\setminus_{\{\tau_{i}\}_{I}})\rightarrow
\mathbb{R}^{k+q}$ and write the impulse condition (\ref{b:1}) with
boundary condition (\ref{b:2}) in the next form
$
\mathfrak{L} x(\cdot,\varepsilon)= \delta,
$
where $\delta:=[\begin{array}{c}
                  \gamma \\
                  \alpha
                \end{array}]\in \mathbb{R}^{k+q}.$

Thus we got the weakly perturbed boundary value problem for
integro-differential system instead weakly perturbed impulsive
boundary problem (\ref{b:1}), (\ref{b:2}):
\begin{equation}\label{b:7}
  \dot{x}(t)-\Phi(t)\int\limits_{a}^b
\Big[A(s)x(s)+B(s)\dot{x}(s)\Big]ds=f(t)+\int\limits_{a}^b
K(t,s)ds\cdot u,\
\end{equation}

\begin{equation}\label{b:8}
\mathfrak{L} x(\cdot,\varepsilon)= \delta \in \mathbb{R}^{k+p},
\end{equation}
$$t\in [a,b]\setminus_{\{\tau_{i}\}_{I}},\tau_{i}\in
(a,b),~~i=1,...,p.$$

The corresponding generating problem $(\varepsilon=0)$ is unsolvable

\begin{equation}\label{b:9}
  \dot{x}(t)-\Phi(t)\int\limits_{a}^b
\Big[A(s)x(s)+B(s)\dot{x}(s)\Big]ds=f(t),\
\end{equation}

\begin{equation}\label{b:10}
\mathfrak{L} x(\cdot,\varepsilon)= \delta\in \mathbb{R}^{k+q}.
\end{equation}

Then according to \cite{BoiHol}, we can formulate the following
criterion for the solvability of boundary-value problem (\ref{b:9}),
(\ref{b:10}).

\begin{theorem}
 {Let  ${\rm rank}Q=n_{2}\leq min(k+q,r_{1})$. The
homogeneous boundary-value problem (\ref{b:9}),
(\ref{b:10})~($f(t)=0,\delta=0$) possesses exactly $r_2$
($r_{2}=r_{1}-n_{2}$) linearly independent solutions of the form:

$$x(t,c_{r_{2}})=\Psi_0(t)P_{D_{r_{1}}}P_{Q_{r_{2}}}c_{r_{2}},~
c_{r_{2}}\in \mathbb{R}^{r_{2}},$$

$$r_{1}=m+n-rank D,~~~r_2=m+n-rank D-rank Q.$$
The inhomogeneous problem (\ref{b:9}), (\ref{b:10})  is solvable if
and only if $f(t)\in L_2[a,b]$ and $\delta\in \mathbb{R}^{k+q}$
satisfy conditions:
\begin{equation}
\label{b:11}P_{D_{d_{1}}^{*}}
\tilde{b}=0,~P_{Q_{d_{2}}^{*}}(\delta-\mathfrak{L}(F(\cdot)))=0,
\end{equation}

$$d_{1}=m-rank D,~d_{2}=k+q-rank Q.$$
In this case, the problem (\ref{b:9}), (\ref{b:10})  possesses an
$r_{2}$-parameter family of solutions:

$$x(t)=\Psi_0(t)P_{D_{r_{1}}}P_{Q_{r_{2}}}c_{r_{2}}+\Psi_{0}(t)P_{D_{r_{1}}}Q^{+}
(\delta-\mathfrak{L}(F(\cdot)))+F(t),$$ where
$Q=\mathfrak{L}X_{r_1}(\cdot)$~is an $(k+q)\times r_1$ matrix,~the
matrix~$Q^{+}$ is pseudoinverse (in the Moore–Penrose sense, \cite{BoiSam})
to the matrix
Q,~$F(t)=\tilde{f}(t)+\Psi_0(t)D^{+}\tilde{b},$~$X_{r_{1}}(t)=\Psi_0(t)P_{D_{r_{1}}}$
is an $n\times r_1$
matrix,~$D=\Big[I_{m}-\int\limits_{a}^b[A(s)\Psi(s)+B(s)\Phi(s)]ds,
-\int\limits_{a}^b A(s)ds\Big]$ is an $m\times(m+n)$
matrix.}\end{theorem}

Here, $\Psi(t)= \int\limits_{a}^t \Phi(s) ds$,
 $\Psi_0(t)=\Big[\Psi(t),I_{n}\Big]$,~$\tilde{b}=\int\limits_{a}^b[A(s)\tilde{f}(s)+B(s)f(s)]ds.$
\linebreak $P_{D}, P_{D^{*}}$ are are $(m+n)\times (m+n),~m \times
m~$matrices (orthoprojectors)  projecting $\mathbb{R}^{m+n}$ and
$\mathbb{R}^{m}$ onto $N(D)=\ker D$ and onto $N(D^{*})=$ker
$D^{*}=$coker $D$, respectively, that is, $P_{D}: \mathbb{R}^{m+n}
\rightarrow N(D),$ $P^{2}_{D}=P_{D}=P^{*}_{D},$ and $P_{D^{*}}:
\mathbb{R}^{m} \rightarrow N(D^{*}),$
$P^{2}_{D^{*}}=P_{D^{*}}=P^{*}_{D^{*}}.$ The matrix $P_{D_{r_1}}
(P_{D^{*}_{d_1}})$ is composed of a complete system of $r_1$ ($d_1$)
linearly independent columns (rows) of the matrix $P_{D}
(P_{D^{*}}).$ $P_{Q}, P_{Q^{*}}$ are $r_1\times r_1,~(k+q) \times
(k+q)~$matrices (orthoprojectors)  projecting $\mathbb{R}^{r_1}$ and
$\mathbb{R}^{k+q}$ onto $N(Q)=\ker Q$ and onto $N(Q^{*})=$ker
$Q^{*}=$coker $Q$, respectively, that is, $P_{Q}: \mathbb{R}^{r_1}
\rightarrow N(Q),$ $P^{2}_{Q}=P_{Q}=P^{*}_{Q},$ and $P_{Q^{*}}:
\mathbb{R}^{k+q} \rightarrow N(Q^{*}),$
$P^{2}_{Q^{*}}=P_{Q^{*}}=P^{*}_{Q^{*}}.$ The matrix $P_{Q_{r_2}}$
($P_{Q^{*}_{d_2}}$) is composed of a complete system of $r_2$
($d_2$) linearly independent columns (rows) of the matrix $P_{Q}
(P_{Q^{*}}).$

Consider the case when one of the conditions ($\ref{b:11}$) is not
fulfilled. Then the boundary-value problem
($\ref{b:9}$),~($\ref{b:10}$) does not have solutions.

It is of interest to analyze whether it is possible to make problem
($\ref{b:9}$),~($\ref{b:10}$) solvable by introducing linear
perturbation and (in the case of positive answer to this \linebreak
question) determine perturbation $K(t,s)$ and $K_1(t,s)$ required to
make the boundary-value problem ($\ref{b:7}$),~($\ref{b:8}$)
everywhere solvable.

Using the solvability criterion ($\ref{b:11}$) of the linear inhomogeneous boundary value problem, we obtain the solvability condition for problem (\ref{b:7}), (\ref{b:8})
\begin{equation}
\label{bondar:912}
P_{D^{*}_{d_1}}\tilde{b_1}=0,\,\,P_{Q^{*}_{d_2}}\big\{\delta-\mathfrak{L} F_1(\cdot)\big\}=0,
\end{equation}
$$d_1=m-rank D,\,\, d_2=p-rank Q.$$
We know
\begin{equation*}
\tilde{b_1}=\tilde{b}+ \big( \int\limits_{a}^b
\Big[A(s)\int\limits_{a}^s \int\limits_{a}^b K(\tau,s)dsd\tau+B(s)\int\limits_{a}^b K(s,\tau)d\tau \Big]ds\big ) u
\end{equation*}
\begin{multline*}
F_1(t)=F(t)+\int\limits_{a}^t \int\limits_{a}^b K(t,s)ds dt+\Psi_{0}(t)D^{+}\int\limits_{a}^b
\Big[A(s)\int\limits_{a}^s \int\limits_{a}^b K(\tau,s)dsd\tau+\\
+B(s)\int\limits_{a}^b K(s,\tau)d\tau \Big]ds\cdot u,
\end{multline*}
get the following algebraic system for definition $u:$
\begin{equation}\label{bondar:913}
\big( P_{D^{*}_{d_1}}\int\limits_{a}^b
\Big[A(s)\int\limits_{a}^s \int\limits_{a}^b K(\tau,s)dsd\tau+B(s)\int\limits_{a}^b K(s,\tau)d\tau \Big]ds\big) u=-P_{D^{*}_{d_1}}\tilde{b},
\end{equation}
\begin{multline}\label{bondar:914}
\big(P_{Q^{*}_{d_2}}\ell\int\limits_{a}^\cdot \int\limits_{a}^b K(t,s)ds dt+\Psi_{0}(\cdot)D^{+}\int\limits_{a}^b
\Big[A(s)\int\limits_{a}^s \int\limits_{a}^b K(\tau,s)dsd\tau+\\
+B(s)\int\limits_{a}^b K(s,\tau)d\tau \Big]ds\big) u=P_{Q^{*}_{d_2}}\big\{\delta-\mathfrak{L} F(\cdot)\big\}.
\end{multline}
Here we have $P_{D}, P_{D^{*}}$~--- $(m+n)\times (m+n)$ and $m \times
m$-dimensional matrices, orthoprojectors that acting from $R^{m+n}$ and $R^{m}$
to  to the kernel and cokernel of matrix $D,$ respectively. The matrix $P_{D_{r_1}}
(P_{D^{*}_{d_1}})$ consists of the complete system $r_1$ ($d_1$)
linearly of independent columns (rows) of the matrix $P_{D} (P_{D^{*}});$
the matrix $Q=\ell X_{r_1}(\cdot)$~--- $p\times r_1$ is dimensional,~$Q^{+}$ is pseudo-inverse in the Moore-Penrose sense of the matrix Q \cite{BoiSam}.
$P_{Q}, P_{Q^{*}}$~--- $r_1\times r_1$ and $p \times
p$-dimensional matrices, orthoprojectors acting from $R^r_1$ and $R^p$
to the kernel and co-kernel of the matrix $Q,$ respectively. The matrix $P_{Q_{r_2}}
(P_{Q^{*}_{d_2}})$ consists of the complete system $r_2$ ($d_2$)
f linearly independent columns (rows) of the matrix $P_{Q} (P_{Q^{*}}).$

Then combining (\ref{bondar:913}), (\ref{bondar:914}) we obtain the following system
\begin{equation}\label{bondar:915}
Uu=g,
\end{equation}
where
$(d_1+d_2)\times n$~---he dimensional matrix $U$ has the form
\begin{equation}\label{bondar:916}
U:=\left[%
\begin{array}{c}
  P_{D^{*}_{d_1}}\int\limits_{a}^b
\Big[A(s)\int\limits_{a}^s \int\limits_{a}^b K(\tau,s)dsd\tau+B(s)\int\limits_{a}^b K(s,\tau)d\tau \Big]ds\\
\\
  P_{Q^{*}_{d_2}}\ell\int\limits_{a}^\cdot \int\limits_{a}^b K(t,s)dsdt+\Psi_{0}(\cdot)D^{+}\int\limits_{a}^b
\Big[A(s)\int\limits_{a}^s \int\limits_{a}^b K(\tau,s)dsd\tau+\\
+B(s)\int\limits_{a}^b K(s,\tau)d\tau \Big]ds\\
\end{array}%
\right]
\end{equation}
$(d_1+d_2)\times 1$ the dimensional vector $g$ is given as follows
\begin{equation}\label{bondar:917}
g:=\left[%
\begin{array}{c}
  -P_{D^{*}_{d_1}}\tilde{b}\\
   P_{Q^{*}_{d_2}}\big\{\delta-\mathfrak{L} F(\cdot)\big\}\\
\end{array}%
\right],
\end{equation}
System (\ref{bondar:915}) is solvable if and only if the condition is
 \begin{equation}\label{bondar:918}
P_{U^{*}}g=0
\end{equation}
and has a solution $ u=U^{+}g+P_{U}c,\,\,c\in \mathbb{R}^{n}.$ Here $U^{+}$ is pseudo-inverse (according to Moore–Penrose) to $U$ is $n\times (d_1+d_2)$ dimensional matrix  $P_{S^{*}}$ is $(d_1+d_2)\times (d_1+d_2)$ dimensional matrix (orthoprojector) that projects $\mathbb{R}^{d_1+d_2}$ onto $N(S^{*}),$ $P_{S}$ is an $(n\times n)$ dimensional matrix (orthoprojector) that projects $\mathbb{R}^{n}$ onto $N(S).$ The following theorem is true.

\begin{theorem}\label{bon_t1} Impulse system of integro-differential equations \eqref{b:1}, \eqref{b:2}, which is unsolvable for $u=0$ and for
$\forall f(t)\in L_{2}[a,b]$ will have a solution if and only if the following condition
	$$P_{U^{*}}g=0$$
is fulfilled.
	In this case, the control variable  $u$ should be selected as follows:
	$$u=U^{+}g+P_{U}c,   c \in R^{n}.$$		
\end{theorem}

\textbf{Remark.} Under condition  $P_{U^{*}}g=0$ the control of $u \in R^{n}$ may not be unique, because it depends on an arbitrary constant
 $P_{U}c \in R^{n}.$ This makes it possible to use this control to investigate problems that are often encountered in the theory of optimal control.



\begin{flushleft}
\textbf{3. Application}
\end{flushleft}
In economics, impulse systems of integro-differential equations with control variables can be applied to model and analyze various economic phenomena. One specific example is modeling economic policy interventions in a macroeconomic context. Let's consider a simplified scenario.\\
\textbf{Example: Economic Policy Intervention}\\
Suppose we have a simplified macroeconomic model that describes the dynamics of an economy.
The state variable $x(t)$ represents the level of economic output, and we want to study the impact of a government's fiscal policy intervention on output.

The inhomogeneous system of integro-differential equations with impulsive actions can be formulated as follows:
\begin{equation*}
	\dot{x}(t)-\Phi(t)\int\limits_{a}^b
	\Big[A(s)x(s)+B(s)\dot{x}(s)\Big]ds=f(t)+\int\limits_{a}^b K(t,s)ds u,~~t\epsilon [a,b],
\end{equation*}
In this context:\\
-- $x(t)$ represents the economic output at time $t;$\\
-- $\dot{x}(t)$ represents the rate of change of economic output;\\
-- $\Phi(t)$ represents a matrix describing how economic output depends on various factors;\\
-- $A(s)$ and $B(s)$ matrices capture the historical dependencies of economic output;\\
-- $f(t)$ represents external economic factors, such as government spending or taxation;\\
-- $K(t, s)$ is a kernel function that models the impulse effect at specific times;\\
-- $u$ represents a control variable that can be adjusted by policymakers.
Now, let's discuss a specific application related to economic policy.

\textbf{Scenario: Government Stimulus Package}\\
Suppose a government wants to stimulate economic growth during a recession. They decide to implement a stimulus package that includes increased government spending $(f(t))$ and a tax cut. The control variable $u$ represents the magnitude of the tax cut, and policymakers want to determine the optimal value of $u$ to maximize economic output over a certain time period.\\
\textbf{1. Model Calibration:} Researchers can calibrate the parameters of the model, including $\Phi(t),$ $A(s),$ $B(s),$ and the kernel function $K(t, s),$ using historical economic data and econometric techniques. This step ensures that the model accurately captures the dynamics of the economy.\\
\textbf{2. Objective Function:} Define an objective function that represents the government's goal, such as maximizing economic output over a specified time horizon. This objective function depends on the control variable $u.$\\
\textbf{3. Optimization:} Use optimization techniques, such as dynamic programming or numerical optimization, to find the optimal value of $u$ that maximizes the objective function. This step involves solving the integro-differential equations for different values of $u$ and selecting the one that yields the highest economic output.\\
\textbf{4. Policy Implementation:} Implement the recommended tax cut $(u)$ as part of the stimulus package. Monitor the economic performance over time.\\
\textbf{5. Evaluation:} Continuously evaluate the impact of the policy intervention on economic output. Adjust the tax cut $(u)$ as needed based on real-time economic data and feedback.

This example illustrates how impulse systems of integro-differential equations with control variables can be applied to model and optimize economic policy interventions. The control variable $(u)$ represents the policy instrument that policymakers can adjust to achieve specific economic goals, such as stimulating growth during a recession. Researchers and policymakers can use this modeling approach to make data-driven decisions and assess the effectiveness of economic policies.

In summary, research on impulse systems of integro-differential equations offers a wide range of possibilities for both theoretical advancements and practical applications. It is an interdisciplinary field with potential contributions to mathematics, engineering, and various other scientific disciplines. As technology continues to advance, the understanding and control of complex dynamic systems, including those with impulsive behaviors, will remain an important area of study.

\emph{\textbf{The publication contains the research results of project No.~~2020.02/0089 with the grant support of the National Research Fund of Ukraine.
}}

Ivanna Bondar

Institute of Mathematics

The National Academy of Science of Ukraine

Tereshchenkivs'ka str. 3

01 601 Kyiv

UKRAINE

E-mail: bondar.ivnn@gmail.com

\end{document}